\documentclass[sn-mathphys-num]{sn-jnl}

\usepackage{graphicx}%
\usepackage{multirow}%
\usepackage{amsmath,amssymb,amsfonts}%
\usepackage{amsthm}%
\usepackage{mathrsfs}%
\usepackage[title]{appendix}%
\usepackage{xcolor}%
\usepackage{textcomp}%
\usepackage{manyfoot}%
\usepackage{booktabs}%
\usepackage{algorithm}%
\usepackage{algorithmicx}%
\usepackage{algpseudocode}%
\usepackage{listings}%

\usepackage{graphicx}
\usepackage{amsfonts}
\usepackage{amsmath}
\usepackage{amssymb}
\usepackage{mathtools}
\usepackage{mathabx}
\usepackage{csquotes}

\theoremstyle{thmstyleone}%
\theoremstyle{thmstyletwo}%

\theoremstyle{thmstylethree}%

\begin{document}

\title[A lower bound for the modulus of the Dirichlet eta function ..]{A lower bound for the modulus of the Dirichlet eta function  on a partition $\mathcal{P}$ from 2-D principal component analysis and transitive composition}


\author[1.2]{\fnm{Author} \sur{Yuri Heymann}}\email{y.heymann@yahoo.com}



\affil[1]{ \emph{Address in Switzerland:} 9 rue Chantepoulet, 1201 Geneva}
\affil[2]{\orgdiv{School of Mathematics}, \orgname{Georgia Institute of Technology}, \orgaddress{\street{686 Cherry Street}, \city{Atlanta}, \postcode{30332-0160}, \state{GA}, \country{USA}}}



\abstract{The present manuscript aims to derive an expression for the lower bound of the modulus of the Dirichlet eta function on vertical lines $\Re(s)=\alpha$.  The approach employs concepts of two-dimensional principal component analysis built on a parametric ellipse, to match the dimensionality of the complex plane. The one-sided lower bound $\forall s \in \mathbb{C}$ s.t. $\Re(s) \in \mathcal{P}$, $| \eta(s) | \geq \left| 1 - \frac{\sqrt{2}}{2^\alpha} \right|$, where $\eta$ is the Dirichlet eta function, is related with the Riemann hypothesis as $|\eta(s)| > 0$ for any $s \in \mathbb{C}$ s.t. $\Re(s) \in \mathcal{P}$, where $\mathcal{P}$ is a partition spanning one half of the critical strip depending upon a variable. We propose the composite lower bound $\forall s \in \, \mathbb{C}$ s.t. $\Re(s) \in \,]1/2,1[$,  $|\eta(s)| \geq \text{Min}\left(1- \frac{\sqrt{2}}{2^{\alpha}},\frac{\sqrt{2}}{2^\alpha}-\frac{\sqrt{2}}{2}\right)$, resulting from transitive composition in $\eta(s) = \left(1-\frac{2}{2^s} \right) \zeta(s)$. As a founding principle, the solution space of the set of solutions referring to such $\mathcal{L}^2$-problem is a representation of the space spanned by explanatory variables satisfying its algebraic form.}




\maketitle

\section{Introduction}

The Dirichlet eta function is an alternating series related to the Riemann zeta function of interest in the field of number theory for the study of the distribution of primes \cite{Riemann1859}. Both series are tied together on a two-by-two relationship expressed as $\eta(s) = \left(1 - 2^{1-s} \right) \, \zeta(s)$ where $s$ is a complex number. The location of the non-trivial zeros of the Riemann zeta function in the critical strip $\Re(s) \in \, ]0,1[$ is key in prime-number theory. For example, the Riemann-von Mangoldt explicit formula, as an asymptotic expansion of the prime-counting function, involves a sum over the non-trivial zeros of the Riemann zeta function \cite{Mangoldt}. Riemann hypothesis, the scope of the domain of existence of the zeros in the critical strip has implications for the accurate estimate of the error involved in the prime-number theorem and a variety of conjectures such as the Lindel\"{o}f hypothesis \cite{Edwards}, conjectures about short intervals containing primes \cite{Dudek}, Montgomery's pair correlation conjecture \cite{Montgomery}, the inverse spectral problem for fractal strings \cite{Lapidus}, etc. Moreover, variants of the Riemann hypothesis falling under the generalized  Riemann hypothesis in the study of modular L-functions \cite{Milinovich} are core for many fundamental results in number theory and related fields such as the theory of computational complexity. For instance, the asymptotic behavior of the number of primes less than $x$ described in the prime-number theorem, $\pi(x) \sim \frac{x}{\ln(x)}$, provides a smooth transition of time complexity as $x$ approaches infinity. As such, the time complexity of the prime- counting function using the $x/(\ln x)$ approximation is of order $\mathcal{O}(M(n) \log n)$, where $n$ is the number of digits of $x$ and $M(n)$ is the time complexity for multiplying two n-digit numbers. This figure is based on the time complexity to compute the natural logarithm with the arithmetic-geometric mean approach where $n$ represents the number of digits of precision.

The below definitions are provided on an informal basis as a supplement to standard definitions when referring to reals, complex numbers and holomorphic functions. A complex number is the composite of a real and imaginary number, forming a 2-D surface, where the purely imaginary axis is represented by the letter $i$ such that $i^2=-1$. The neutral element and index $i$ form a basis spanning some sort of vector space. 
The Dirichlet eta function is a holomorphic function having for its domain a subset of the complex plane where reals are positive and denoted $\mathbb{C}^{+}$, which arguments are sent to codomain in $\mathbb{C}$. As a holomorphic function it is characterized by its modulus, a variable having for support an axial vertex orthonormal to the complex plane. 
The conformal way to describe space in geometrical terms is the orthogonal system, consisting of eight windows delimited by the axes of the Cartesian coordinates, resulting from the union of the three even surfaces of Euclidean space. 

As an extension of the former, a Hilbert space is a multidimensional space which, in the current context of some functions in $\mathcal{L}^2$ space referring to squared-integrable functions, is comprised of wave functions expressing components or basis elements of space, further equipped with an inner product defining a norm and angles between these functions. The well-known Riemann hypothesis, part of the eighth problem David Hilbert presented at the International Congress of Mathematicians in Paris in 1900 \cite{Hilbert}, has a strong relation with such Hilbertian spaces, as unveiled in the remaining of the manuscript.

\vspace{1mm}

\rule{23mm}{0.3pt}

\vspace{1mm}

\noindent
As a reminder, the definition of the Riemann zeta function and its analytic continuation to the critical strip are displayed below. The Riemann zeta function is commonly expressed as follows:

\begin{equation}
\zeta(s) = \sum_{n=1}^{\infty} \frac{1}{n^{s}} \,,
\end{equation}

\noindent
where $s$ is a complex number and $\Re(s)>1$ by convergence in the above expression. 

\vspace{3mm}
\noindent
The standard approach for the analytic continuation of the Riemann zeta function to the critical strip $\Re(s) \in \, ]0, 1[$ is performed with the multiplication of $\zeta(s)$ with the function $\left(1-\frac{2}{2^s} \right)$, leading to the Dirichlet eta function.  By definition, we have:

\begin{equation}
\eta(s) = \left(1-\frac{2}{2^s} \right) \zeta(s) =  \sum_{n=1}^{\infty} \frac{(-1)^{n+1}}{n^{s}} \,,
\end{equation}

\noindent
where $\Re(s)>0$ and $\eta$ is the Dirichlet eta function. By continuity as $s$ approaches one, $\eta(1)=\ln(2)$. 

\vspace{3mm}

\noindent
The function $\left(1-\frac{2}{2^s} \right)$ has an infinity of zeros on the line $\Re(s)=1$ given by $s_k=1+\frac{2 k \pi i}{\ln 2}$ where $k \in \, \mathbb{Z}^{*}$. As $\left(1-\frac{2}{2^s} \right) = 2 \times \left(2^{\alpha-1} \, e^{i \, \beta \ln 2} -1 \right) / \,  {2^\alpha \, e^{i \, \beta \ln 2}}$, the factor $\left(1-\frac{2}{2^s} \right)$  has no poles nor zeros in the critical strip $\Re(s) \in \, ]0, 1[$. As such, the Dirichlet eta function can be used as a proxy of the Riemann zeta function for zero finding in the critical strip $\Re(s)\in \, ]0, 1[$. 

\vspace{2mm}

\noindent
From the above, the Dirichlet eta function is expressed as:

\begin{equation}
\eta(s) =  \sum_{n=1}^{\infty} \frac{ (-1)^{n+1} \, e^{-i \, \beta \ln(n)}}{n^\alpha} \,,
\end{equation}

\noindent
where $s= \alpha + i \, \beta$ is a complex number, $\alpha$ and $\beta$ are real numbers. 

\vspace{2mm}
\noindent
 We have $\frac{1}{n^{s}} = \frac{1}{n^\alpha \exp(\beta i \ln n)} =  \frac{1}{n^\alpha \left(\cos(\beta \ln n) + i \sin(\beta \ln n) \right)}$. We then multiply both the numerator and denominator by $\cos(\beta \ln(n))-i \sin(\beta \ln(n))$. After several simplifications, $\eta(s) =  \sum_{n=1}^{\infty}  \frac{(-1)^{n+1} \left[\cos(\beta \ln n)-i \sin(\beta \ln n) \right]}{n^\alpha}$.

\vspace{3mm}
\noindent
In the remainder of the manuscript, the Riemann zeta function refers to its formal definition and analytic continuation by congruence. 

\section{Mathematical development}
\label{s2}

\subsection{Elementary propositions no. 1 $\sim$ 4} 
\label{s2.1}

\hfill\break

\vspace{1mm}
\noindent{\textbf{Proposition 1}} Given $z_1$ and $z_2$ two complex numbers, we have:

\begin{equation}
\Big|z_1 + z_2 \,\Big| \geq \Big| |z_1| - |z_2|\,\,\Big| \,,
\end{equation}

\noindent
where $|z|$ denotes the modulus of the complex number $z$. This  is the well-known reverse triangle inequality, which is valid for any normed vector space (including complex numbers), where the norm is subadditive over its domain of definition see \cite{Dragomir2004,NakaiTada1995}.

\vspace{2mm}

\noindent{\textbf{Proposition 2}} Let us consider an ellipse $(x_t,y_t) = \left[ a \cos(t), b \sin (t) \right]$ where $a$ and $b$ are two positive reals corresponding to the lengths of the semi-major and semi-minor axes of the ellipse ($a \geq b$) and $t \in [0, 2 \pi]$ is a variable having a correspondence with the angle between the x-axis and the vector $(x_t,y_t)$.

\vspace{2mm}
\noindent
Let us set $t$ such that the semi-major axis of the ellipse is aligned with the x-axis, which is the angle maximizing the objective function defined as the modulus of $(x_t,y_t)$. When $|(x_t,y_t)|$ is maximized, we have:

\begin{equation}
|(x_t,y_t)|=x_t+y_t = a \,.
\end{equation}

\noindent
Note that by maximizing $x_t+y_t$, we would get $|(x_t,y_t)| < x_t+y_t$, as the expression $x_t+y_t$ is maximized when $t=\arctan(b/a)$, leading to $\text{max}(x_t+y_t)=\sqrt{a^2+b^2}$. 

\vspace{2mm}
\noindent
\textit{Background and geometry.} Principal Component Analysis (PCA), is a statistical concept for reducing the dimensionality of a variable space by representing it with a few orthogonal variables capturing most of the variability of an observable. An ellipse centered on the origin of the coordinate system can be parametrised as follows: $(x_t,y_t) = \left[ a \cos(t), b \sin (t) \right]$ where $a$ and $b$ are positive real numbers corresponding to the lengths of the semi-major and semi-minor axes of the ellipse ($a \geq b$) and $t \in \mathbb{R}$ is a variable having correspondence with the angle between the x-axis and the vector $(x_t,y_t)$. The objective function $|(x_t,y_t)|$ is maximized with respect to $t$ when the major axis is aligned with the x-axis. When $|(x_t,y_t)|$ is at its maximum value, we have $|(x_t,y_t)|=x_t+y_t=a$ where $a$ is the length of the semi-major axis. The modulus of $(x_t,y_t)$ is as follows: $|(x_t,y_t)| = \sqrt{a^2 \cos^2(t) + b^2 \sin^2(t)} \leq a, \forall t \in \mathbb{R}$. In the principal component analysis representation, $x_t$ represents the first principal component and $y_t$ the second principal component. Suppose $x_t$ and $y_t$ were not orthogonal, then there would be a non-zero phase shift $\varphi$ between the components, i.e. $x_t= a\cos(t)$ and $y_t=b \sin(t+\varphi)$. 

\vspace{2mm}

\noindent{\textbf{Proposition 3}} Given a vector $V_{\mathcal{E}}=[u(\mathcal{E}), v(\mathcal{E})]$ defined in a bidimensional vector space, where $u(\mathcal{E})$ and $v(\mathcal{E})$ are two real-valued functions say on $\mathbb{R}^{\mathcal{\nu}} \rightarrow \mathbb{R}$, where $\nu$ represents the degrees of freedom of the system. The reference of a point in such system is described by the set $\mathcal{E}=\{\varepsilon_1, \varepsilon_2,..., \varepsilon_{\nu}\}$ representing a multidimensional coordinate system. Thus, we have:

\begin{equation}
|V_{\mathcal{E}}|=u(\mathcal{E})+v(\mathcal{E}) \,,
\end{equation}

\noindent

if and only if $u(\mathcal{E}) \, v(\mathcal{E}) = 0$ and $u(\mathcal{E})+v(\mathcal{E}) \geq 0$. 

\vspace{2mm}

\noindent
Given a basis set $\{e_1 , e_2\}$ of the above-mentioned vector space, where $|e_i|=1$ for $i=1, 2$, \textit{proposition 3} is true if and only if the inner product across basis elements is equal to zero, i.e. $e_1 \cdot e_2 = 0$, meaning that the basis elements are disentangled from each other. We say that $\{e_1 , e_2\}$ is an orthonormal basis. This condition is also necessary for $\textit{propositions 2} \, and \,\textit{4}$ to be true, in a 2-D Cartesian frame. 

\begin{proof}
\noindent
By the square rule, we have $\left( u+ v \right)^2 = u^2+v^2 + 2 \, u \,v$. The modulus of a vector $V$ as defined in such a two-dimensional frame is $|V|=\sqrt{u^2+v^2}$, leading to $|V|=|u+v|$ if and only if $u \, v = 0$, which is provided as a complement to \textit{proposition 2}. The above as a support of pre-Hilbertian spaces by the scalar product $u \, v$ is a prerequisite for Hilbertian spaces of squared-integrable functions referring to such $\mathcal{L}^2$ spaces equipped of an inner product.
\end{proof}

\vspace{2mm}

\noindent{\textbf{Proposition 4}} Given a circle of radius $r \in \, \mathbb{R}^+$ parametrized as follows: $(x_t,y_t)=[r \, \cos(t), r \, \sin(t)]$ where $t$ is a real variable in $[0, 2 \pi]$, we construct a function $f(t)=a \cos(t) + b \sin (t+\varphi)$ where $a$ and $b$ are two positive reals and $\varphi$ is a real variable which can be positive or negative such that:

\begin{equation}
r \, \cos(t) + r \, \sin(t) = a \cos(t) + b \sin (t+\varphi) \,,
\end{equation}

\noindent
$\forall t \in \mathbb{R}$ and where $\varphi$ is a real variable of $t$ (when $a$ and $b$ are scalars). 

\vspace{2mm}

\rule{23mm}{0.3pt}

\vspace{3mm}

\noindent
As an excerpt of below proof elements $\forall t \in \, [0, 2 \pi]$ and $\forall \delta \in \, [-r, r]$ we have: 
\vspace{2mm}
\\ 
$r\cos(t)+r \sin(t)= (r+\delta) \cos(t) + \sqrt{r^2+\delta^2} \sin(t+ \varphi)$, where $\varphi=-\arctan{\delta/r}$, making $\varphi$ independent of $t$ by some correspondence between $a$, $b$ and $\varphi$ represented as parametric functions of $\delta$ and where $\delta/r \in [-1,1]$.

\vspace{2mm}
\noindent
As such $a=r+\delta$ and $b=\sqrt{r^2+\delta^2}$ where $\delta=- r \, \tan(\varphi)$, yielding $a/r=1-\tan \varphi$ and $b/r=\sqrt{1+\tan^2 \varphi}$ where $\varphi \in \, ]-\pi/4,\pi/4 [$.

\vspace{2mm}

\rule{46mm}{0.3pt} 

\vspace{3mm}

\noindent
Say $u_t=a \cos(t)$ and $v_t=b \sin (t+\varphi)$. 

\vspace{2mm}
\noindent
When $a \geq b$, the first component $u_t$ is the one carrying most of the variance\footnote{The variance in the statistical meaning as applied to a $\mathbb{R} \times \mathbb{R}$ pure function $f: t \mapsto f(t)$ represented as a univariate $f_t=f(t)$ is measured by an estimator, e.g. $\text{Var}(f(t))=\frac{1}{|b-a|}\int_{a}^{b} \, f^2(t)\, dt$ over interval $[a,b]$. The variance of a component of $f(t)$, is defined by its attribution of variance with respect to total variance.} of expression $f(t)$, meaning it is the leading component. Thus, we have:
\begin{equation}
|(x_t, y_t)|  \leq \text{max}(v_t) \leq \text{max}(u_t) \,,
\end{equation}

\noindent
$\forall t \in \,[0, 2 \pi]$, where $\text{max}(u_t)$ is the maximum value of $u_t$ and $\text{max}(v_t)$ the maximum value of $v_t$ over the interval $[0,2\pi]$.

\vspace{2mm}

\noindent
When $a \leq b$, the component $v_t$ carries most of the variance of $f(t)$, meaning it is the leading component, and we have:

\begin{equation}
 \text{max}(u_t) \leq |(x_t, y_t)|  \leq \text{max}(v_t) \,,
\end{equation}

\noindent
$\forall t \in \,[0, 2 \pi]$, where $\text{max}(u_t)$ and $\text{max}(v_t)$ as above.

\vspace{3mm}

\noindent
When the functions $u_t$ and $v_t$ are orthogonal, i.e. $\varphi=0$, we have $r=a=b$.

\begin{proof}
Given that $(r+\delta) \cos(t) + (r-\delta_2) \sin(t) = r \cos(t) + r \sin(t) + \delta \cos(t) - \delta_2 \sin(t)$, where $\delta$ and $\delta_2$ are variables in $\mathbb{R}$ i.e. sensitive to $t$. As we want  $\delta \cos(t) - \delta_2 \sin(t)=0$, we have $\delta_2= \delta \cot(t)$. Thus, we get: $r \cos(t) + r \sin(t) = (r+\delta) \cos(t) + \left( r - \delta \cot(t)\right) \sin(t)$.  As $r \sin(t) -\delta \cos(t) = \sqrt{r^2+\delta^2} \sin(t+\varphi)$ where $\varphi=-\arctan{\delta/r}$, we get $\forall t \in \, [0, 2 \pi]$, $r\cos(t)+r \sin(t)= (r+\delta) \cos(t) + \sqrt{r^2+\delta^2} \sin(t+ \varphi)$. We set $a=r+\delta$ and $b=\sqrt{r^2+\delta^2}$ where $-r \leq \delta \leq r$, leading to (7). If $\delta \geq 0$, we have $r \leq \sqrt{r^2+\delta^2} \leq r+\delta$, leading to (8). If $\delta \leq 0$,  we have $r+\delta \leq r \leq \sqrt{r^2+\delta^2}$, leading to (9). As $\text{max}\{r \cos(t) + r \sin(t)\}= \sqrt{2} \, r$ which occurs when $t= \frac{\pi}{4}$, we have $\delta \in [-r, r]$ for any $r \geq 0$.
\end{proof}

\subsection{The Riemann hypothesis, background and implications}

According to the RH (Riemann hypothesis), all non-trivial zeros lie on the critical line $\Re(s)=1/2$, meaning that RH is true if and only if $|\eta(s) |> 0$ for any $\alpha \in \, ]0,1/2[$ and $]1/2,1[$. The one-sided RH test is a consequence of the Riemann zeta functional, leading to \textit{prop 7} in \cite{Heymann}, i.e. Given $s$ a complex number and $\bar{s}$ its conjugate, if $s$ is a zero of the Riemann zeta function in the strip $\Re(s) \in \, ]0, \, 1[$, then $1-\bar{s}$ is also a zero of the function. 

The converse is also true, meaning if $s$ is not a zero of the Riemann zeta function in the critical strip, then $1-\bar{s}$ is not a zero as well. The one-sided test is enough for the RH to be true. When $\alpha = 1/2$, $|\eta(s)| \geq 0$ means the Dirichlet eta function can have some zeros on the critical line $\Re(s) = 1/2$, which is known to be true \cite{Titchmarsh}, p. 256. As the Dirichlet eta function and the Riemann zeta function share the same zeros in the critical strip, we have to show that $|\eta(s) |> 0$ for any $s$ in the critical strip not on the critical line to say the Riemann hypothesis is true. 

In the present study, we propose several lower bounds for the modulus of the Dirichlet eta function that are one-sided lower bounds (i.e. apply on one side of the critical strip). These are eqns. (17), (18) and (20).  Note the zeros of the Dirichlet eta function on the line $\Re(s) = 1$, produce violations of (17) in the neighborhood of such zeros. Lower bounds (18) and (20) suggest that $|\eta(s)|>0$ when $\alpha \in ]1/2,1[$ as per the one-sided RH test. The asymmetrical skew behind the Riemann hypothesis is examined in section 4, by analysing the solution set of a pair of Taylor polynomials.


\section{The lower bound of the Dirichlet eta modulus as a floor function}
\label{s3}

In standard notations, the point $s$ is expressed as $s= \alpha + i \, \beta$ where $\alpha$ and $\beta$ are real in their corresponding basis belonging to $\mathbb{C}$. 

\vspace{2mm}
\noindent
Note the zeros of the Dirichlet eta function and its complex conjugate are the same. For convenience, we introduce the conjugate of the Dirichlet eta function, expressed as follows:

\begin{equation}
\widebar{\eta}(s) = \sum\limits_{n=1}^{\infty} (-1)^{n+1} \frac{e^{i \, \beta \ln n}}{n^\alpha} \,,
\end{equation}

\noindent
where $\Re(s) >0$.  By applying the reverse triangle inequality to (10) (see \textit{proposition 1}), we get:

\begin{equation}
\begin{split}
\left|\widebar{\eta}(s)\right| & \geq  \left| 1- \left| \sum\limits_{n=2}^{\infty} (-1)^{n+1} \frac{e^{i \, \beta \ln n}}{n^\alpha} \right| \right| \,,
\end{split}
\end{equation}

\noindent
where $|z|$ denotes the norm of the complex number $z$. 

\vspace{3mm}
\noindent
With respect to expression $\left| \sum\limits_{n=2}^{\infty} (-1)^{n+1} \frac{e^{i \, \beta \ln n}}{n^\alpha} \right|$, its decomposition into sub-compon- ents $u_n=\frac{(-1)^{n+1}}{n^\alpha} e^{i \beta \ln n}$ is a vector representation where $\beta \ln n + (n+1) \pi$ is the angle between the real axis and the orientation of the vector itself, and where $\frac{1}{n^\alpha}$ is its modulus. The idea is to apply a rotation by an angle $\theta$ to all component vectors simultaneously, resulting in a rotation of the vector of their sum. The resulting vector after rotation $\theta$ expressed in Euler's notation is $v_{\theta, \beta}=\sum\limits_{n=2}^{\infty} \frac{(-1)^{n+1}}{n^\alpha} e^{i (\, \beta \ln n+ \theta)}$, where $\theta$, $\alpha$ and $\beta$ are real numbers. 

\vspace{3mm}
\noindent
Let us introduce the objective function $w$, defined as the sum of the real and imaginary parts of $v_{\theta, \beta}$, i.e. $w = v_x+ v_y$ where $v_x = \Re(v_{\theta, \beta})$ and $v_y= \Im(v_{\theta, \beta})$. We get:

\begin{equation}
\begin{split}
w= & \sum\limits_{n=2}^{\infty} \frac{(-1)^{n+1}}{n^\alpha} \left( \cos(\beta \ln n + \theta) + \sin(\beta \ln n + \theta) \right) \\
&= \sum\limits_{n=2}^{\infty} (-1)^{n+1} \frac{\sqrt{2}}{n^\alpha} \cos\left(\beta \ln n + \theta - \frac{\pi}{4}\right) \,.
\end{split}
\end{equation}

\noindent
The trigonometric identity $\cos(x)+\sin(x)=\sqrt{2}\cos\left(x- \frac{\pi}{4} \right)$ which follows from $\cos(a) \cos(b) + \sin(a) \sin(b) = \cos(a-b)$ with $b=\frac{\pi}{4}$ is invoked in (12), see \cite{Stegun} formulas 4.3.31 and 4.3.32, p. 72. The finite sum of a subset of the elements of the second line of (12) from 2 to $n \in \mathbb{N}$ is further referred to as the $w$-series.

\vspace{3mm}
\noindent
While orthogonality between vectors is defined in terms of the scalar product between such pairs, for real functions on $\mathbb{R}$ to $\mathbb{R}$ we usually define an integration product forming an $\mathcal{L}^2$ space. Let us say we have two real-valued functions $f$ and $g$, which are squared-$Lebesgue$ integrable on a segment $[a,b]$ and where the inner product between $f$ and $g$ is given by:

 \begin{equation}
\langle f,g \rangle = \int_{a}^{b} f(x) \, g(x) dx  \,.
\end{equation}

\noindent
The functions $f$ and $g$ are squared-$Lebesgue$ integrable, meaning such functions can be normalized; i.e. the squared norm as defined by $\langle f, f \rangle$ is finite. For sinusoidal functions such as sine and cosine, it is common to say $[a, b] = [0, 2 \,\pi]$, which interval corresponds to one period. The condition for functions $f$ and $g$ to be orthogonal is that the inner product as defined in (13) is equal to zero.

\vspace{3mm}
\noindent
We proceed with the decomposition of the objective function $w$ into dual components $w_1$ and $w_2$, expressed as follows:

\vspace{1mm}

\begin{equation}
w_1= - \frac{\sqrt{2}}{2^\alpha} \cos\left(\beta \ln(2) + \theta - \frac{\pi}{4}\right) \,,
\end{equation}
\noindent
and
\begin{equation}
w_2 =  \sum\limits_{n=3}^{\infty} (-1)^{n+1} \frac{\sqrt{2}}{n^\alpha} \cos\left(\beta \ln n + \theta - \frac{\pi}{4}\right) \,,
\end{equation}

\noindent
where most of the variance of the $w$-series comes from the leading component in a direction of $\mathcal{L}^2$ space.

\vspace{3mm}

\noindent
By construction $w$ is the sum of the real and imaginary parts of $v_{\theta, \beta}$ which are orthogonal functions. Let us say $v_{\theta}$ is the parametric notation of $v_{\theta, \beta}$ for a given $\beta$ value, and $\alpha$ implicitly. We note that for any given $\beta$ value, the complex number $v_{\theta}$ describes a circle in the complex plane, centered on the origin. Hence, in symbolic notations, $w$ can be written as $w=r \cos(t) + r \sin(t)$. The components $w_1$ and $w_2$ can be expressed as $w_1=a \cos(t)$ and $w_2=b \sin(t+\varphi)$ where $t$ is a variable in $[0, 2 \pi]$ and $\varphi$ some variable in $\mathbb{R}$. As we suppose that $w_1$ carries most of the variance of $w$ (i.e. $a \geq b$),  the modulus $|v_{\theta}|$ is smaller or equal to the maximum value of $|w_1|$, by \textit{proposition 4}. Yet,  $|v_{\theta}|$ is equal to $\text{max}\{w_1\}$, if  $w_1$ and $w_2$ are orthogonal and $t=0$.  By \textit{proposition 2}, at its maximum value $|(w_1,w_2)| =w_1+w_2$ when $\varphi = 0$ and $t=0$, which in light of the above, is also equal to the maximum value of $|v_{\theta}|$. As the inner product between $w_1$ and $w_2$ does not depend on $\theta$, orthogonality between $w_1$ and $w_2$ is determined by $\beta$ values. We then apply a rotation by an angle $\theta$ to maximize the objective function $|(w_1,w_2)|$. We consider two complementary scenarios respectively, depending on whether $w_1$ is the leading component of the $w$-series or some other function; i.e. $w_2$ as the alternative, or in a direction of $\mathcal{L}^2$ space.


\vspace{4mm}
\noindent{\textbf{When $w_1$ is the leading component:}} 
\vspace{2mm}

\noindent
Suppose $w_1$ is the leading component in some regions of the critical strip denoted $\mathcal{A}$. In this scenario, as we suppose $w_1$ and $w_2$ are orthogonal at the maximum value of $|(w_1,w_2)|$, i.e. $w_1=\frac{\sqrt{2}}{2^\alpha}$ and $w_2=0$, we get $\textit{max}\{|(w_1,w_2)|\} = \frac{\sqrt{2}}{2^\alpha}$, also equal to $\text{max}\left|v_{\theta,\beta}\right|$ by \textit{proposition 2}. \\
If we suppose that $w_1$ and $w_2$ are not orthogonal and form an acute angle, by \textit{proposition 4} eq. (8) we would get $\left| \sum_{n=2}^{\infty} (-1)^{n+1}\frac{e^{i \, \beta \ln n}}{n^\alpha}\right| < \frac{\sqrt{2}}{2^\alpha}$, and $| \eta(s) |$ would be strictly larger than zero when $\alpha=1/2$ in (11). This would imply that the Dirichlet eta function does not have zeros on the critical line $\Re(s)=1/2$, which is known to be false. Hence, we can say that when $w_1$  is the leading component, the functions $w_1$ and $w_2$ are orthogonal  at the maximum value of $|v_{\theta}|$, i.e. $w_1$ acts as first principal component. 
As we suppose that:

\begin{equation}
 \forall s \in \, \mathbb{C^{+}} \text{   s.t. } s \, \in \mathcal{A},\, \left| \sum_{n=2}^{\infty} (-1)^{n+1}\frac{e^{i \, \beta \ln n}}{n^\alpha}\right| \leq \frac{\sqrt{2}}{2^\alpha} \,,
\end{equation}

\noindent
where $\mathbb{C^{+}}$ is the subset of $\mathbb{C}$ s.t. $\Re(s)>0$ (acronym s.t. standing for "such that"). 

\vspace{2mm}
\noindent
Say for $\Re(s)=\alpha \geq \frac{1}{2}$, (11) and (16) imply that :

\begin{equation}
| \eta(s) | \geq \left| 1 - \frac{\sqrt{2}}{2^\alpha} \right|\,,
\end{equation}

\noindent
for any $s \in \mathcal{A}$ s.t. $\Re(s)=\alpha \in \, [1/2,\infty[$, where $\mathcal{A} \subseteq \mathbb{C^{+}}$. This is a one-sided lower bound function $L_\alpha$ which is always true when (16) applies.

\noindent
By the conjugate, we have:

\begin{equation}
| \eta(1-s) | \geq \left| 1 - \frac{\sqrt{2}}{2^{1-\alpha}} \right|\,,
\end{equation}

\noindent
which is proposed as an alternative one-sided lower bound applicable when (17) fails, i.e. as the left-side of the critical strip for $Re(s) \leq 1/2$ is not impacted by the zeros of the Dirichlet eta function on the line $\Re(s)=1$.

\vspace{2mm}
\noindent
As such lower bound $\forall s \in \mathbb{C}$ s.t. $\Re(s) \in \mathcal{P}$, $| \eta(s) | \geq \left| 1 - \frac{\sqrt{2}}{2^\alpha} \right|$, where $\eta$ is the Dirichlet eta function is related with the Riemann hypothesis as $|\eta(s)| > 0$ for any $s \in \mathbb{C}$ s.t. $\Re(s) \in \mathcal{P}$, where $\mathcal{P}$ is a partition spanning one half of the critical strip on either side of the critical line $\Re(s)=1/2$ depending upon a variable. 

\vspace{4mm}

\noindent{\textbf{When $w_1$ is not the leading component, i.e. $w_2$ as the alternative, or in a direction of $\mathcal{L}^2$ space:}} 
\vspace{2mm}

\noindent
As a complementary of the former, say $\mathcal{A}^c$ is some regions of the critical strip where $w_1$ is not leading. In this scenario, we fall on eq. (9) of \textit{proposition 4}, yielding: 

\begin{equation}
 \forall s \in \, \mathbb{C^{+}} \text{ s.t. } s \, \in \mathcal{A}^c, \, \left| \sum_{n=2}^{\infty} (-1)^{n+1}\frac{e^{i \, \beta \ln n}}{n^\alpha}\right| > \frac{\sqrt{2}}{2^\alpha} \,,
\end{equation}

\noindent
where $\mathbb{C^{+}}$ is defined as above and $w_2$ is the strictly leading component, hence (17) is no longer guaranteed to work and we need some additional thinking. 

\begin{figure}[h]
 \includegraphics[width=6.5cm, height=6.5cm]{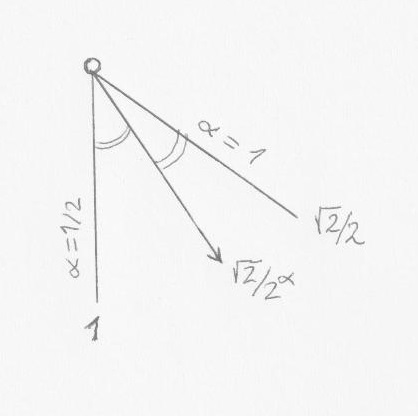}
  \caption{Pendulum model representing the gaps when $\alpha$ approaches $1$ on the right-hand side of the critical strip $\Re(s) \geq 1/2$. } 
\end{figure}

\vspace{2mm}
\noindent
The function $\left(1-\frac{2}{2^s} \right)$ has an infinity of zeros on the line $\Re(s)=1$ given by $s_k=1+\frac{2 k \pi i}{\ln 2}$ where $k \in \, \mathbb{Z}^{*}$. The Dirichlet eta function as a product of the former and the Riemann zeta function has some zeros on the line $\Re(s)=1$. A violation of (17) by a tiny epsilon was reported by Vincent Granville in the neighborhood of a point given by $s=0.75+580.13 \, i$. This point is located near a zero of the Dirichlet eta function on the line $\Re(s)=1$, for $k=64$. From the pendulum model (see Fig. 1), cases when $\alpha$ approaches $1$ producing a violation of (17) in the neighborhood of zeros on the line $\Re(s)=1$ are taken care of by replacing the wider gap $1-\frac{\sqrt{2}}{2^\alpha}$ by the more narrow gap $\frac{\sqrt{2}}{2^\alpha}-\frac{\sqrt{2}}{2}$. The resulting composite lower bound is

\begin{equation}
\left| \eta(s)\right| \geq \text{Min}\left( 1-\frac{\sqrt{2}}{2^\alpha} , \frac{\sqrt{2}}{2^\alpha} - \frac{\sqrt{2}}{2}\right) \,,
\end{equation}

\noindent
where $\alpha \in ]1/2,1[$, from transitive composition in $\eta(s) = \left(1-\frac{2}{2^s} \right) \zeta(s)$. Case when $\alpha = 1$ shows no positive gap (i.e. $|\eta(s)| \geq 0$) due to the zeros on the line $\Re(s)=1$. 

\vspace{2mm}
\noindent
The narrow gap function in the pendulum model is an application of the reverse triangle inequality to function $\frac{1}{\sqrt{2}} \times \left(1 - \frac{2}{2^s} \right)$. The validity of this gap function follows from the observation that $\forall k \in \, \mathbb{Z^{*}}$, when $(\beta \, \ln 2) \in \,\, ]2 k \pi, 2k\pi +2\pi[$, $\exists \,\delta >0$ such that $\forall \alpha \in \, ]1-\delta,1]$, $|\zeta(s)| \, \sqrt{\left( 1 - \frac{2}{2^\alpha} \, \cos(\beta \, \ln 2) \right)^2 +\left( \frac{2}{2^\alpha} \, \sin \left(\beta \, \ln 2\right)\right)^2} \geq \frac{1}{\sqrt{2}} \, \left| 1 - \frac{2}{2^{\alpha}}\right|$. From the Laurent series of the Riemann zeta function on the line $\Re(s)=1$, we get $\zeta(1+i\, \beta)=-\frac{1}{\beta} \, i + \sum_{n=0}^{\infty} \frac{\gamma_n}{n !} (i \, \beta)^n$, where $\gamma_n$ are the Stieltjes constants. With an expansion up to the first five Stieltjes constants, $|\zeta(1+i\, \beta)|$ is minimized around $\beta \simeq 2.9$ and equal to $0.8419 \pm 0.001$ at the minimum. Hence, for all $\beta \in \mathbb{R}$, we can say that $|\zeta(1+i\, \beta)| > \frac{1}{\sqrt{2}}$.  This motivates the application of the reverse triangle inequality to $\frac{1}{\sqrt{2}} \times \left(1 - \frac{2}{2^s} \right)$ resulting from $\eta(s) = \left(1-\frac{2}{2^s} \right) \zeta(s)$. (20) results from transitive composition of the Dirichlet eta function. A consequence of (20) is that its validity implies the Riemann hypothesis is true, by first principle.

\section{The roots of a biquadratic form and higher orders - the asymmetrical skew behind the Riemann hypothesis}
\label{s4}

While the Riemann zeta function has no zeros when $\Re(s)>1$ (see \cite{Ivic}) and no zeros on the line $\Re(s) = 1$ (a result of Hadamard and de la Vall{\'e}e Poussin proofs of the prime-number theorem \cite{Poussin,Hadamard}), remains to be proven that all zeros in the critical strip $\Re(s) \in \, ]0,1[$ lie on the critical line $\Re(s)=1/2$ as from the Riemann hypothesis.

From the Riemann zeta functional $\zeta(1-s) = \Gamma(s) (2 \pi)^{s-1} 2 \cos \left( \frac{s \pi}{2} \right) \zeta(s) $ by Riemann (details of the derivation in \cite{Edwards}, p.13.), we obtain a theorem:  ``Given $s$ a complex number and $\bar{s}$ its complex conjugate, if $s$ is a zero of the Riemann zeta function in the strip $\Re{s} \in \, ]0,1[$, then $1-\bar{s}$ is also a zero.'' (see prop. 8 in \cite{Heymann} for detailed derivation).

We use the Dirichlet eta function as a proxy of the Riemann zeta function for zero finding in the critical strip. Many attemps to prove the Riemann hypothesis with the Riemann zeta functional have failed due to the zero divided by zero singularity at zeros of the function. We find that the ``0/0 problem'' (zero divided by zero singularity) with the Riemann zeta functional, is resolved with a Taylor expansion of $\varphi_n(\alpha)=\frac{1}{n^{\alpha}}$ and $\phi_n(\alpha)=\frac{1}{n^{1-\alpha}}$ around $1/2+\Delta\alpha$ (allowing the construction of polynomials about the axis $\Re(s)=1/2$).

\vspace{2mm}
\noindent
Given a weight function $w_n \in [-1,1]$ where $n=1,...,N$ and $N \in \, \mathbb{N}$, the question of whether the number of non-trivial roots such that $\sum_{n=2}^{\infty} \frac{w_n}{n^\alpha}=\sum_{n=2}^{\infty} \frac{w_n}{n^{1-\alpha}}=\lambda_i$ is a countable set, where trivial roots lie on the critical line $\alpha=1/2$ arises.

\vspace{2mm}
\noindent
By Taylor expansion of $\varphi_n(\alpha)=\frac{1}{n^{\alpha}}$ in $\alpha_0=1/2$+$\Delta \alpha$, we have:

\begin{equation}
\varphi_n(\alpha) = \varphi_n(1/2)+\varphi_n^{\prime}(1/2) \, \Delta \alpha + \frac{\varphi_n^{\prime\prime}(1/2)}{2!} \, \Delta \alpha^2 +.... \,.
\end{equation}

\vspace{2mm}
\noindent
By Taylor expansion of $\phi_n(\alpha)=n^{\alpha-1}$ in $\alpha_0=1/2$+$\Delta \alpha$, we have:

\begin{equation}
\begin{split}
\phi_n(\alpha) & = \phi_n(1/2)+\varphi_n^{\prime}(1/2) \, \Delta \alpha + \frac{\phi_n^{\prime\prime}(1/2)}{2!} \, \Delta \alpha^2 +.... \\
& =  \varphi_n(1/2) - \varphi_n^{\prime}(1/2) \, \Delta \alpha + \frac{\varphi_n^{\prime\prime}(1/2)}{2!} \, \Delta \alpha^2 +.... \,.
\end{split}
\end{equation}

\noindent
From (21) and (22), $\sum_{n=2}^{\infty} \frac{w_n}{n^\alpha}=\sum_{n=2}^{\infty} \frac{w_n}{n^{1-\alpha}}$ is rewritten as follows:

\begin{equation}
\sum_{n=2}^{\infty} w_n \, \varphi_n^{\prime}(1/2)+ \sum_{n=2}^{\infty} w_n \, \frac{\varphi_n^{(3)}(1/2)}{3!} \, \Delta \alpha^2 + \sum_{n=2}^{\infty} w_n \, \frac{\varphi_n^{(5)}(1/2)}{5!} \, \Delta \alpha^4 +... =0 \,,
\end{equation}

\noindent
which expansion is convergent when $|\Delta \alpha| < 1$ (convergent within the critical strip).

\vspace{2mm}
\noindent
As a first guess, we suppose the higher order terms of the Taylor expansion are negligible, meaning (23) is asymptotic to a biquadratic equation $Q(x)=a_4 \, x^4+a_2\, x^2 + a_0$. Let the auxiliary variable be $z = \Delta \alpha^2$, leading to $Q(z)=a_4 \, z^2+a_2\, z + a_0$. For any pair of distinct roots, there is at most one real root of the biquadratic having a value equal to $\lambda=1$ (see Fig. 2). As such the uniqueness of the roots in (23) is given by the discriminant $\Delta=0$ in the biquadratic scenario. Unicity of the roots does not necessarily hold for higher orders. 

\begin{figure}[h]
 \includegraphics[width=12cm, height=6.5cm]{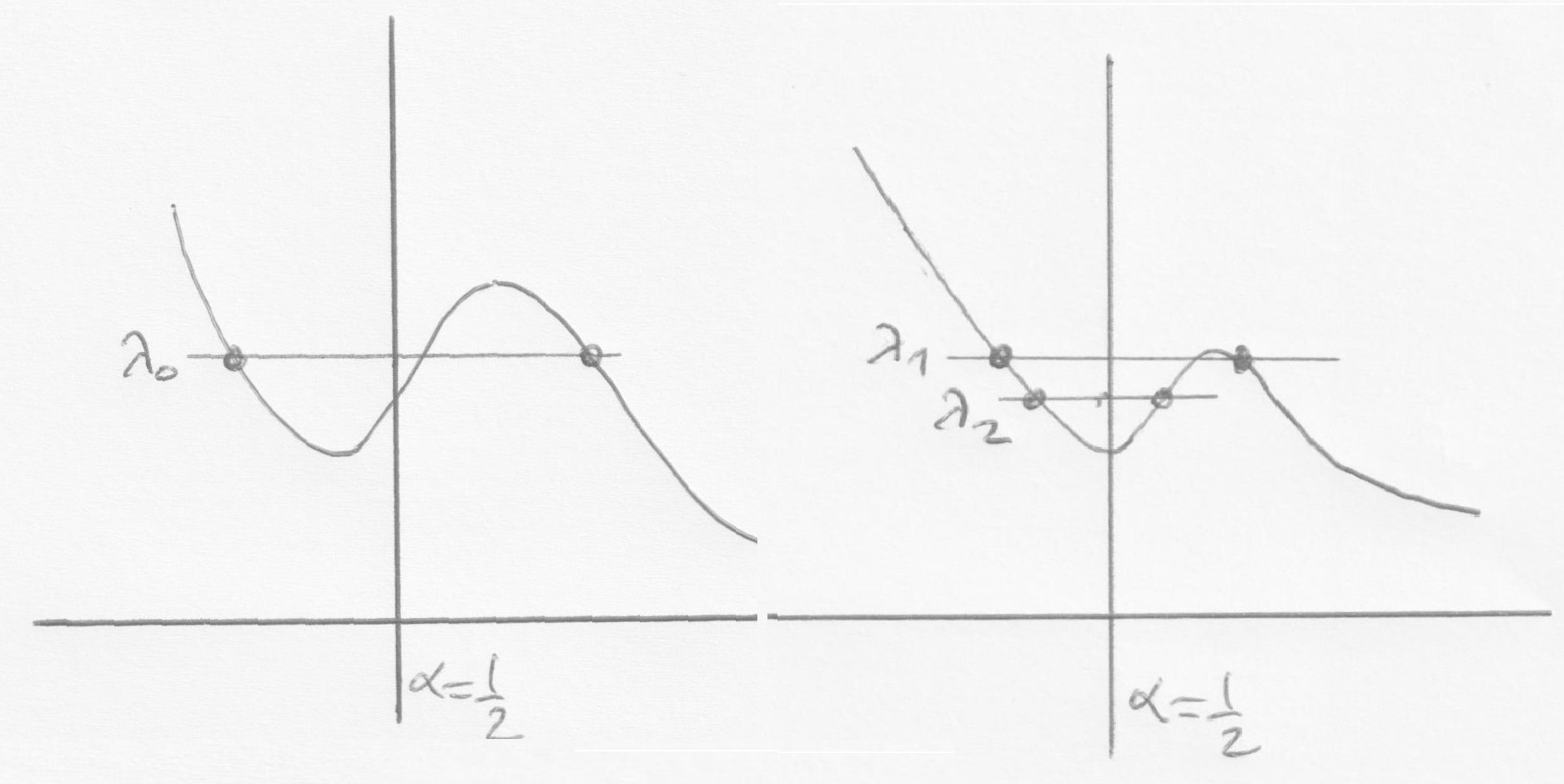}
  \caption{Roots of the biquadratic equation Q(x)=0. By the Riemann zeta functional, the roots of this biquadratic are the set of points equidistant to axis $\Re(s)=1/2$. The graph on the left-hand side shows a simple root, while the graph on the right two distinct roots. The case with two roots is characterized by distinct heights $\lambda_1$ and $\lambda_2$.} 
\end{figure}

\vspace{2mm}
\noindent
Any point in the critical strip beside the trivial case $\alpha=1/2$ is a zero of the Dirichlet eta function, provided the two conditions set forth below are satisfied simultaneously (see theorem \textit{prop 8} in \cite{Heymann}, from Riemann zeta functional):

\begin{equation}
\sum_{n=2}^{\infty} (-1)^{n}\, \frac{\cos\left( \beta\, \ln n\right)}{n^\alpha}= \sum_{n=2}^{\infty} (-1)^{n}\, \frac{\cos\left( \beta\, \ln n\right)}{n^{1-\alpha}} = 1 \,,
\end{equation}

\noindent
and

\begin{equation}
\sum_{n=2}^{\infty} (-1)^{n}\, \frac{\sin\left( \beta\, \ln n\right)}{n^\alpha}= \sum_{n=2}^{\infty} (-1)^{n}\, \frac{\sin\left( \beta\, \ln n\right)}{n^{1-\alpha}} = 0 \,.
\end{equation}

\noindent
Stemming out of \textit{prop 8}, these two equalities stand for the real and imaginary parts of the Dirichlet eta function to be equal to zero simultaneously, respectively on the right-hand side of the critical strip, i.e. $\alpha \in ]1/2,1[$ and conjugated form $1-\alpha$ (on the left). This entails that both polynomials relative to (24) and (25) as obtained by Taylor expansion of $\varphi_n(\alpha)=\frac{1}{n^{\alpha}}$ and $\phi_n(\alpha)=n^{\alpha-1}$ in $\alpha_0=1/2$ have to be satisfied simultaneously at a value of $\Delta \alpha$ for a given point to be a zero of the Dirichlet eta function, whereas both sets of weights in (24) and (25)  are orthogonal with each other by the sine-cosine orthogonality.

\vspace{2mm}
\noindent
The factors relative to both polynomials as seen in (23) are expressions of the form $\sum_{n=2}^{\infty} a_n \, \cos(\theta_n)$ and $\sum_{n=2}^{\infty} a_n \, \sin(\theta_n)$ respectively, where $a_n$ and $\theta_n$ are real series indexed by $n$ natural as functions of $\beta$. Hence, we can say there exists reals $L$ and $\Psi$, such that $\sum_{n=2}^{\infty}a_n \cos(\theta_n)=L \, \cos(\Psi)$ and $\sum_{n=2}^{\infty}a_n \sin(\theta_n)=L \, \sin(\Psi)$. The former system of biquadratic equations extended by Taylor expansions up to order $n$, leads to its canonical form expressed as follows:

\begin{equation}
L_0 \, \cos(\Psi_0) + L_1 \, \cos(\Psi_1) \, z + L_2 \cos(\Psi_2) \, z^2 +... + L_n \, \cos(\Psi_n) z^n= 0 \,,
\end{equation}

\noindent
and

\begin{equation}
L_0\, \sin(\Psi_0) +L_1\, \sin(\Psi_1) \, z + L_2\, \sin(\Psi_2) \, z^2 +... +  L_n \, \sin(\Psi_n) z^n = 0 \,,
\end{equation}

\noindent
where $z=\Delta \alpha^2$, and $L_j$ and $\Psi_j$ are sequences representing real variables (some functions of $\beta$), where index $j=0,..,n$ and $n \in \mathbb{N}$. Series in (26) and (27) are convergent when $|z| <1$, meaning convergent within the critical strip.

\vspace{2mm}
\noindent
We can easily show that the system of polynomials (26) and (27) admits for solution points on the critical line as given by $z=0$ and $L_0 =0$ and a set of unfeasible solutions. By a linear combination of (26) and (27), we get that $\forall \varphi \in \, \mathbb{R}$:

\begin{equation}
L_0 \, \sin(\Psi_0+\varphi) +L_1\, \sin(\Psi_1+\varphi) \, z + L_2\, \sin(\Psi_2+\varphi) \, z^2 +... +  L_n \, \sin(\Psi_n+\varphi) z^n = 0 \,,
\end{equation}

\noindent
which expansion is convergent when $|z| <1$ (within the critical strip).

\begin{proof}
From the identity $\sin(a+b)=\sin(a) \cos(b)+\cos (a) \, \sin(b)$, we get that 
\begin{equation*}
\begin{split}
L_0 \, \sin(\Psi_0+\varphi) +L_1\, \sin(\Psi_1+\varphi) \, z + L_2\, \sin(\Psi_2+\varphi) \, z^2 +... = \\
\cos \varphi \times \left(L_0 \, \sin(\Psi_0) +L_1\, \sin(\Psi_1) \, z + L_2\, \sin(\Psi_2) \, z^2 +... \right) + \\
+\sin \varphi  \times \left(L_0 \, \cos(\Psi_0) +L_1\, \cos(\Psi_1) \, z + L_2\, \cos(\Psi_2) \, z^2 +... \right) = 0
\end{split}
\end{equation*}
which is true  $\forall \varphi \in \, \mathbb{R}$.
\end{proof}

\noindent
By excluding the case $z=0$ and $L_0=0$ (zeros on the critical line), a trivial solution satisfying (28) is when all Latin letters $L_0,L_1,...,L_n$ in (26, 27) are equal to zero simultaneously. As an example, the biquadratic scenario where letter $L_0=0$ (expansion up to $z^2$), (26) and (27) as a pair carries two candidate solutions $z=0$ or $z=-L_1/L_2$. We show that such solutions where $z \neq 0$ are not feasible. When $L_0=0$, (28) is satisfied provided $z=0$ (by a linear combination of (26) and (27)). In the case when $L_0 \neq 0$, (28) cannot be satisfied for all $\varphi \in \, \mathbb{R}$ (it is only satisfied for distinct values of $\varphi$ and $z$). The case where all Latin letters $L_0,L_1,...,L_n$ are equal to zero simultaneously, is a hypothetical scenario elaborated further down.

\begin{figure}[h]
 \includegraphics[width=6.5cm, height=5.9cm]{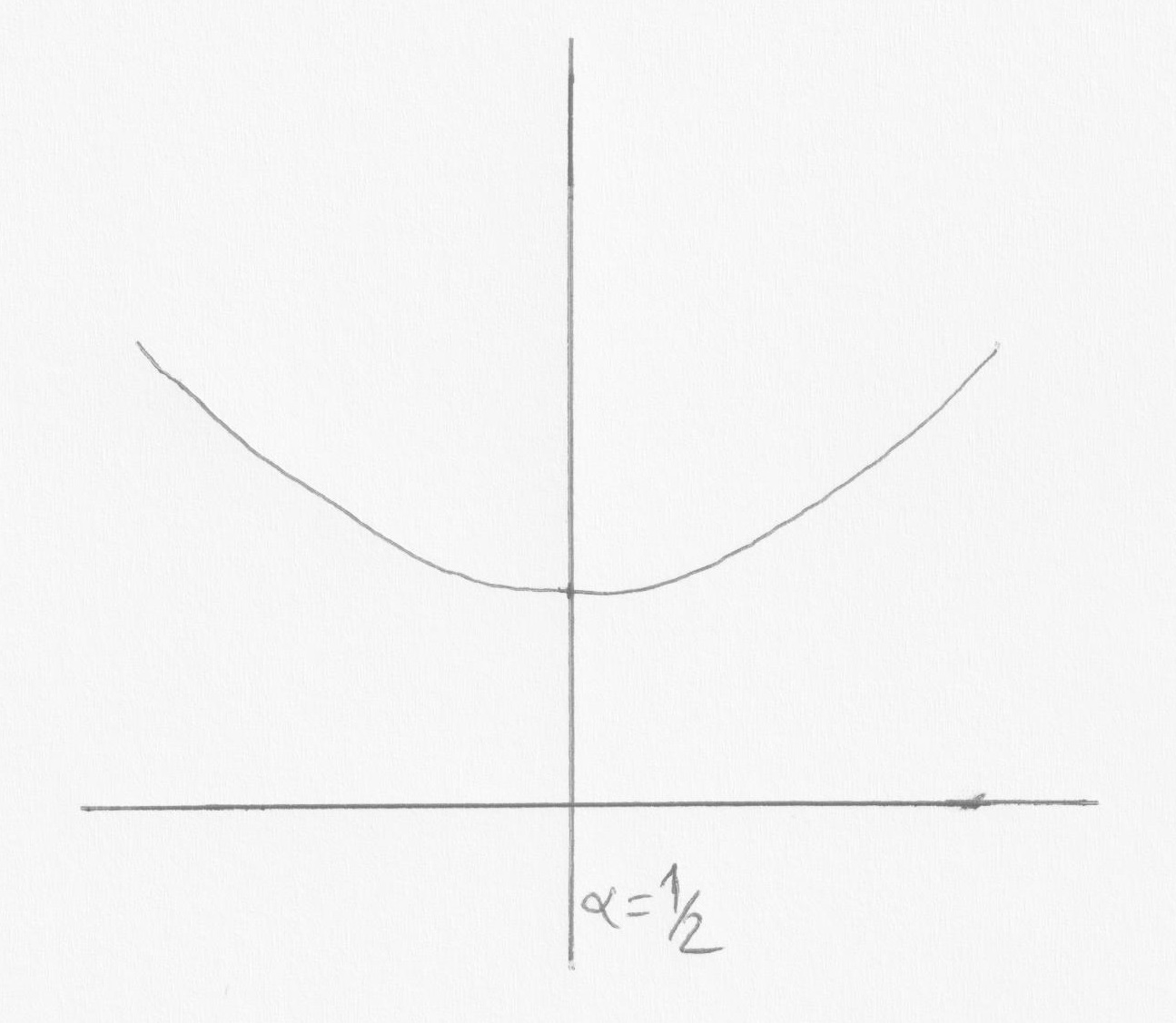}
  \caption{Hypothetical scenario of a symmetrical curve about the vertical axis $\Re(s)=1/2$. As the roots of such biquadratics (extendable to polynomials of arbitrary orders) are given by equidistant points to the axis $\Re(s)=1/2$, the scenario of a symmetrical curve about $\Re(s)=1/2$ yields an infinity of solutions, i.e. any real $\alpha$ satisfies the biquadratic polynomial.} 
\end{figure}

\vspace{2mm}
\noindent
The question arises whether there exists a solution such that (24) or (25) is satisfied for all $\alpha$ in the critical strip (see scenario all Latin letters $L_0,L_1,L_2,..,$etc are equal to zero simultaneously). By the Taylor expansion, the curves spanned by such polynomials are asymptotic to a symmetrical shape about axis $\Re(s)=1/2$, see Fig. 3. This problem can be formulated as follows: Can a linear combination of exponentials $e^{-a_n\, x}$, where $a_n=\ln n$ is a real series indexed by natural $n$, lead to a symmetrical function about the vertical axis $x=1/2$ ?

\vspace{2mm}
\noindent
The asymptotic expression of a linear combination of $e^{-a_n\, x}$ forming a weighted sum, is defined by function $F : \mathbb{R} \to \mathbb{R}$ where $F(x)=\sum_{n=2}^{N} w_n \, e^{-a_n \, x}$, and $w_n$, $n=2,..,N \in \mathbb{N}$ are the weights. This function is symmetrical with respect to vertical axis $x=1/2$ if and only if: 

\begin{equation}
\frac{\partial F(x)}{\partial x} = -\frac{\partial F(1-x)}{\partial x} \,,
\end{equation}

\noindent
for all $x \ge 0$ (as Dirichlet proxy is defined for positive $\alpha$). By expansion into a weighted sum, the former gives $\forall x \in \mathbb{R}^{+}$: %

\begin{equation}
\sum_{n=2}^{N} w_n \, \left(e^{-a_n x}- e^{-a_n\,(1-x)}\right) = 0 \,,
\end{equation}

\noindent
which is only true provided all weights are equal to zero. As such, the construction of a symmetrical function about the vertical axis $x=1/2$ from a linear combination of exponentials is not feasible without having all weights equal to zero.

\vspace{2mm}
\noindent
The zero-weighted sum scenario only occurs when $\beta=0$, for which all weights in (25) are equal to zero. Still, $\left| \sum_{n=2}^{\infty} \frac{(-1)^n}{n^{\alpha}}\right| \leq \frac{1}{2^{\alpha}}<1$ for $\forall \alpha >0$, by upper trailing of the sum of a decreasing alternating series (see \cite{Rainville} p.74), meaning the real part of the Dirichlet eta function is strictly positive when $\Re(s)>0$. 

\section{Discussion}

The present work employs the Dirichlet eta function as a proxy of the Riemann zeta function for zero finding in the critical strip, and describes results about the lower bound of the modulus of the Dirichlet eta function as a floor function. The surface spanned by the  modulus of the Dirichlet eta function is a continuum resulting from the application of a real-valued function over the dimensions of the complex plane, which is a planar representation where the reals form a line continuous to the right intersecting the imaginary axis, and where the square of imaginary numbers are subtracted from zero. As a design aspect, the modulus of the Dirichlet eta function is a holographic function sending a complex number into a real number, whereas the floor function is a projection of the former onto the real axis.

In the common scenario when $w_1$ is the leading component, the floor function of the modulus of the Dirichlet eta function on vertical lines $\Re(s) = \alpha$ does not depend on $\beta$, which is reflected by the linear relationship between $\theta$ and $\beta$ in the cosine argument of the first principal component, as a single term of $w$-series. This is no longer the case, when adding together several terms of the $w$-series as first principal component. As a complementary of the former, scenarios when $w_1$ is not the leading component, i.e. $w_2$ is the alternative, or some direction of $\mathcal{L}^2$ space exists in various parts of the domain. For such scenarios, though there is no straight-forward linear relationship between $\theta$ and $\beta$ of the arguments of component $w_2$ as the alternative, for $L_{\alpha} =\left| 1- \sqrt{2} / 2^{\alpha} \right|$ to qualify as the floor function in all the domain, $L_{\alpha}$ is the floor function of the modulus of the Dirichlet by mirror symmetry with respect to line $\Re(s)=1/2$, c.f. (18) in complement of (17). By the perfect matching principle under mirror symmetry, components $w_1$ and $w_2$ at zeros of the Dirichlet eta function resulting from the continuity to the right of the critical line as given by (18), are meant to match corresponding components at such zeros as given by (16) at the limit to the left of line $\Re(s) = 1/2$. An inversion of the roles of components $w_1$ and $w_2$ when approaching the critical line from both sides, suggests there are no such regions where $w_1$ is not leading, which is contiguous with the critical line.

Special case of a polynomial made up of a subset of the cosine terms of the $w$-series which does not depend on $\beta$ occurs, if there exists such a polynomial that is equal to zero for any $\beta$. For such a polynomial to be first principal component involves subsequent components are also equal to zero, leading to the degenerate case $|v_\theta|=0$. This occurs when $\alpha$ tends to infinity, leading to $w=0$ for $\beta$ real, as a special case of the Dirichlet eta function converging towards unity.

The combination of multiple terms of the $w$-series as principal components involves such components are functions composed of terms of the form \\$a_{n}= \pm \frac{1}{n^\alpha}\cos(\beta \ln(n)+ \theta)$, where $n$ is the index of the corresponding term in $v_{\theta}$. Due to the multiplicity of bivariate collinear arguments in the cosine functions, which comovements are not parallel across the index $n$ (as a finite set), there is no straight-forward bijection between $\theta$ and $\beta$, i.e. a one-degree of freedom relationship, such that all cosine arguments $\beta \ln n + \theta$ of the component are decoupled from $\beta$. As aforementioned, the lower bound of the modulus of the Dirichlet eta function needs to be decoupled from $\beta$, to be a floor function on vertical lines. Moreover, the principal components involved in modulus maximization need to be disentangled for PCA to be applicable, which in the current context means the maximum of $|(w_1, w_2)|$ is reached at orthogonality between $w_1$ and $w_2$, as a prerequisite for $L_{\alpha}$ to be a floor function of the modulus of the Dirichlet eta function. As a rule of thumb, one degree of freedom is needed for every additional principal component, when matching the dimensionality of the variable space in the parametric ellipsoidal model.

While many attempts to prove the Riemann hypothesis with the Riemann zeta functional have failed due to the zero divided by zero singularity at the zeros of the function, in section 4, we propose an approach using a Taylor expansion of $\varphi_n(\alpha)=\frac{1}{n^{\alpha}}$ and $\phi_n(\alpha)=\frac{1}{n^{1-\alpha}}$ around $1/2$, allowing the construction of polynomials about the axis $\Re(s)=1/2$. This system of polynomial admits for solution points on the critical line $\Re(s)=1/2$ and a set of unfeasible solutions resulting from the orthogonality between the weight functions of both polynomials. A solution exists if we can construct a symmetrical function about the vertical axis $x=1/2$ from a linear combination of exponentials (see hypothetical scenario all Latin letters in (26, 27) are equal to zero simultaneously), which is prevented by the positivity of the curve, i.e. all weights are not equal to zero when $\beta \neq 0$. For the case $\beta = 0$, the real part of the Dirichlet eta function never reaches zero.

\section{Conclusion}
\label{s6}

In the present manuscript, we propose a lower bound of the modulus of the Dirichlet eta function coming from transitive composition and concepts of 2-D principal component analysis, which is expressed as $\forall s \in \mathbb{C}$ s.t. $\Re(s) \in \, ]1/2,1[$, $|\eta(s)| \geq \text{Min}\left(1- \frac{\sqrt{2}}{2^{\alpha}},\frac{\sqrt{2}}{2^\alpha}-\frac{\sqrt{2}}{2}\right)$, where $\eta$ is the Dirichlet eta function. As a proxy of the Riemann zeta function for zero finding in the critical strip $\Re(s) \in \, ]0,1[$, the above floor function of the modulus of the Dirichlet eta function provides a venue for the Riemann hypothesis (i.e.that the non-trivial zeros lie on the critical line $\Re(s) = 1/2$). 

Many attempts to prove the Riemann hypothesis with the Riemann zeta functional have failed due to the zero divided by zero singularity at zeros of the function. This problem is resolved here with a Taylor expansion of the Dirichlet representation, i.e. analyzing the solutions of a pair of polynomials obtained by a Taylor expansion of $\varphi_n(\alpha)=\frac{1}{n^{\alpha}}$ and $\phi_n(\alpha)=\frac{1}{n^{1-\alpha}}$ about $\alpha_0=1/2+ \Delta \alpha$, an approach unveiling the asymmetrical skew behind the Riemann hypothesis.

We show that this system admits for solution points on the critical line $\Re(s)=1/2$, and a set of unfeasible solutions resulting from orthogonality between the weight functions relative to both polynomials. We further verify that a symmetrical function about the vertical axis $x=1/2$ cannot be constructed from a linear combination of exponentials (except when $\beta=0$, all weights are zeros), preventing a scenario where all Latin letters of the canonical form are equal to zero simultaneously. 







\bibliography{Zeta}

\end{document}